\documentclass[final,sort&compress,3p,times,12pt]{elsarticle}
\usepackage{amsmath}
\allowdisplaybreaks[4]
\usepackage{amsmath}
\usepackage{amsfonts}
\usepackage{}
\usepackage{amssymb}
\usepackage{dsfont}
\usepackage{amsmath}
\usepackage{graphicx}
\usepackage{subfigure}
\usepackage[justification=centering]{caption}
\usepackage{color}
\usepackage{array}
\usepackage{mathrsfs}
\usepackage{multirow}

\newtheorem{Defi}{Definition}
\newtheorem{example}{Example}
\newtheorem{Theo}{Theorem}
\newtheorem{Assu}{Assumption}

\newtheorem{Lem}{Lemma}

\numberwithin{equation}{section}
 \numberwithin{Lem}{section}
 \numberwithin{Defi}{section}
 \numberwithin{Theo}{section}
 \numberwithin{Pro}{section}
 \numberwithin{Rem}{section}
  \numberwithin{Coro}{section}
  \numberwithin{Fig}{section}

\begin{document}

\begin{frontmatter}

\title{Projected Euler method for stochastic delay differential equation under a global monotonicity condition$^\star$ \tnotetext[1]{This work was supported by National Natural Science Foundation of China (No. 11771163).
}}
\author{Min Li$^a$,\quad Chengming Huang$^{a,b,*}$}
\cortext[cor1]{Corresponding author.\\
\emph{Email addresses}: \texttt{liminmaths@163.com} (M. Li),
\texttt{chengming\_huang@hotmail.com} (C. Huang$^*$ )}.\\
 \address{$^a$School of Mathematics and Statistics, Huazhong University of Science and Technology, Wuhan 430074,
 China\\}
 \address{$^b$Hubei Key Laboratory of Engineering Modeling and Scientific Computing, Huazhong University of Science
 and Technology, Wuhan 430074, China\\}
\date{}
\begin{abstract}
 This paper investigates projected Euler-Maruyama method for stochastic delay differential equations under a global monotonicity condition. This condition admits some equations with highly nonlinear drift and diffusion coefficients.  We appropriately generalized the idea of C-stability and B-consistency given by  Beyn et al. [J. Sci. Comput. 67 (2016), no. 3, 955-987] to the case with delay. Moreover, the method is proved to be convergent with order $\frac{1}{2}$ in a succinct way. Finally, some
 numerical examples are included to illustrate the obtained theoretical results.
\end{abstract}
\begin{keyword}
Stochastic delay differential equation; Projected Euler-Maruyama method; Strong convergence; C-stability; B-consistency
\end{keyword}
\end{frontmatter}

\section{Introduction}
Consider d-dimensional nonlinear stochastic delay differential equations (SDDEs)
\begin{equation}\label{a1}
dX(t)=f(X(t),X(t-\tau))dt+g(X(t),X(t-\tau))dW(t),~~~~t>0,
\end{equation}
with the initial condition given by
\begin{equation}
\{X(\theta):-\tau\leq \theta\leq 0\}=\xi\in C([-\tau,0];\mathbb{R}^{d}).
\end{equation}
Here, $X:[0,T]\times \Omega\rightarrow \mathbb{R}^{d}$ denotes the exact solution to \eqref{a1}, the drift term $f:\mathbb{R}^{d}\times \mathbb{R}^{d}\rightarrow \mathbb{R}^{d}$ and diffusion term $g:\mathbb{R}^{d}\times \mathbb{R}^{d}\rightarrow \mathbb{R}^{d\times m}$. And $W(t):=(W_{1}(t), \cdots, W_{m}(t))^{T}$ is an m-dimensional Wiener process defined on given complete probability space $(\Omega, \mathscr{F},\mathbb{P})$ with a filtration $\{\mathscr{F}_{t}\}_{t\geq 0}$ under usual condition (i.e., it is increasing and right continous, and $\mathscr{F}_{0}$ contains all $\mathbb{P}$-null sets). SDDEs can be seen as a generalization of stochastic differential equations, and they play an important role in many phenomena in physics \cite{Beuter1993Feedback,Eurich1996Noise,Longtin1990Noise}. In terms of well-posedness of the equation, there have been extensive study and application of SDDEs. The well known result is that the global Lipschitz condition and the linear growth condition guarantee the existence and uniqueness of analytical solution (see, \cite{MR1335454,MR2380366}). In 2002, Mao \cite{MR1893197} gave the Khasminskii-type condition for SDDEs where linear growth condition was no longer necessary, and global existence and uniqueness of the solution was proved.

Most of SDDEs can not be solved analytically, so numerical calculation is particularly necessary. In the past two decades, a number of numerical methods were investigated under Lipschitz and linear growth condition (see \cite{MR2075010,MR2276834,MR3400432,MR1799186,MR2211500,MR3345286,MR3281888,MR3669732,MR3429327,MR2921197,MR2772279,MR3131860} and references therein). Limited work has been done in
SDDEs whose coefficients do not satisfy the linear growth condition, and this issue received attention only recently. The mean square stability of $\theta$ methods for SDDEs under
a coupled condition was first studied by Huang \cite{MR3123472}.  In 2018, Guo et. al \cite{MR3803361} considered the truncated Euler-Maruyama method for nonlinear SDDEs under the generalized Khasminskii-type condition, and convergence in $L^{q}$
was also derived. Zhang et al. \cite{MR3758637} established  the convergence of partially truncated Euler-Maruyama method for a class of highly nonlinear SDDEs. All their convergence analyses were under the framework given by Higham et al. \cite{MR1949404}, where complex higher moment estimation  and  continuous time extension of the corresponding numerical scheme should be taken into account.  Recently, Beyn et al. \cite{MR3493491,MR3608332} proposed  projected Euler-Maruyama method, projected Milstein  method for SDEs by studying the  C-stability and B-consistency, which can avoid those processes on the discrete time level. In this way, the convergence analysis can be simplified significantly.

Compared with implicit methods for SDEs, the explicit Euler methods process simpler algebraic
structure, and can reach strong order of convergence $1/2$ with cheaper computational cost. However, Hutzenthaler et. al \cite{MR2795791} proved that strong and weak divergence in finite time of the explicit Euler method for
SDEs with superlinearly growing coefficients. Subsequently, some modified  Euler methods, such as tamed and truncated methods,  were constructed to solve the nonlinear SDEs (see \cite{MR3430145,MR2985171}). The main goal of this paper is to generalize the projected Euler methods for SDDEs with superlinearly growth condition.

An outline of this paper is organized as follows. Some assumptions and projected Euler method are introduced in Section $2$. Section $3$ gives the main convergence theorem under the premise of stochastic C-stability and B-consistency. In Section $4$, C-stability and B-consistency of projected Euler-Maruyama method are studied in detail. In Section $5$, some numerical experiments are carried out to verify the theoretical results.  Finally, some conclusions are drawn in the last section.
\section{Preliminaries}
Most of the notations in this paper come from \cite{MR3493491}. For the sake of simplicity, we let
$$
h=\frac{\tau}{M},~~t_{n}=nh,
$$
and
$$
t_{i-M}=t_{i}-\tau,~~ i=0, 1, \cdots, M,
$$
then there exists a positive integer $N$ such that $t_{N}\leq T, t_{N}+h> T$.
Further, we follow the notation of the space of adapted and square integrable grid functions
$$
\mathcal{G}^{2}(\mathcal{T}_{h}):=\{Z:\mathcal{T}_{h}\times\Omega\rightarrow \mathbb{R}^{d}: Z(t_{n})\in L^{2}(\Omega,\mathscr{F}_{t_{n}}, \mathbb{P}; \mathbb{R}^{d})~ \text{for~ all}~ n=0, 1, \cdots, N\}.
$$
With the help of the preceding notations, we can give the definition of stochastic one-step methods.
\begin{Defi}
For every $t, t+h \in [0, T]$ and $Z\in L^{2}(\Omega,\mathscr{F}_{t},\mathbb{P};\mathbb{R}^{d})$, $\Psi$ satisfies the following measurability and integrability condition:
$$
\Psi(Z,t,h)\in L^{2}(\Omega,\mathscr{F}_{t+h}, \mathbb{P}; \mathbb{R}^{d}),
$$
if
\begin{align*}
&X_{h}(t_{i})=\Psi(X_{h}(t_{i-1}),X_{h}(t_{i-M}),h),~~1\leq i\leq N,\\
&X_{h}(t_{i-M})=\xi(t_{i}-\tau),~~~i=0,\cdots,M,
\end{align*}
then we say grid function $X_{h}\in\mathcal{G}^{2}(\mathcal{T}_{h})$ is yield by the $\text{stochastic one-step method}~ (\Psi, h,\xi)$.
\end{Defi}
Taking one step projected Euler method in \cite{MR3493491} into account, we propose our new projected method $(\Psi^{PEM}, h, \xi)$ for SDDEs \eqref{a1} as follows
\begin{align}\label{a21}
\begin{split}
\overline{X}_{h}^{PEM}(t_{i-1}):=&\min(1,h^{-\alpha}|X_{h}^{PEM}(t_{i-1})|^{-1})X_{h}^{PEM}(t_{i-1}),\\
\overline{X}_{h}^{PEM}(t_{i-M}):=&\min(1,h^{-\alpha}|X_{h}^{PEM}(t_{i-M})|^{-1})X_{h}^{PEM}(t_{i-M}),\\
X_{h}^{PEM}(t_{i}):=&\overline{X}_{h}^{PEM}(t_{i-1})+hf\bigg{(}\overline{X}_{h}^{PEM}(t_{i-1}),\overline{X}_{h}^{PEM}(t_{i-M})\bigg{)}\\
&+g\bigg{(}\overline{X}_{h}^{PEM}(t_{i-1}),\overline{X}_{h}^{PEM}(t_{i-M})\bigg{)}\big{(}W(t_{i})-W(t_{i-1})\big{)},~~1\leq i\leq N,
\end{split}
\end{align}
with $X_{h}^{PEM}(t_{i-M})=\xi(t_{i}-\tau),~i=0, 1, \cdots, M$.

Before proceeding further, let us make the following assumptions.
\begin{Assu}\label{assu1}
There exist
positive constant $L$ and parameter $\eta\in (\frac{1}{2}, \infty)$ such that
\begin{align}\label{a13}
\langle x_{1}-x_{2}, f(x_{1},\bar{x_{1}})-f(x_{2},\bar{x_{2}})\rangle+\eta|g(x_{1},\bar{x_{1}})-g(x_{2},\bar{x_{2}})|^{2}\leq L(|x_{1}-x_{2}|^{2}+|\bar{x_{1}}-\bar{x_{2}}|^{2}).
\end{align}
The above expression is referred as to global monotonicity condition. Moreover, we assume that there is constant $q\in (1,\infty)$ such that
\begin{align}\label{a8}
|f(x,\bar{x})|\vee|g(x,\bar{x})|\leq L(1+|x|^{q}+|\bar{x}|^{q}),
\end{align}
\begin{align}\label{a9}
\begin{split}
&|f(x_{1},\bar{x_{1}})-f(x_{2},\bar{x_{2}})|\vee|g(x_{1},\bar{x_{1}})-g(x_{2},\bar{x_{2}})|\\
&\leq (1+|x_{1}|^{q-1}+|x_{2}|^{q-1}+|\bar{x_{1}}|^{q-1}+|\bar{x_{2}}|^{q-1})\big{(}|x_{1}-x_{2}|+|\bar{x_{1}}-\bar{x_{2}}|\big{)},
\end{split}
\end{align}
for all $x, x_{1}, x_{2}, \bar{x}_{1}, \bar{x}_{2}\in \mathbb{R}^{d}$.
\end{Assu}
\begin{Assu}\label{assu2}
The initial data $\xi$ satisfies
$$
|\xi(u)-\xi(v)|\leq K_{1}|u-v|^{\beta},~~~~-\tau\leq v<u\leq 0,
$$
where $K_{1}>0$ and $\beta\in \big{[}\frac{1}{2},1\big{]}$ are constants.
\end{Assu}
\begin{Assu}\label{assumption1}
For every positive number $R$, there exists a positive constant $K_{R}$ such that
\begin{align*}
|f(x,y)-f(\bar{x},\bar{y})|^{2}\vee|g(x,y)-g(\bar{x},\bar{y})|^{2}\leq K_{R}(|x-\bar{x}|^{2}+|y-\bar{y}|^{2})
\end{align*}
for those $x, y, \bar{x}, \bar{y}\in \mathbb{R}^{d}$ with $|x|\vee |y|\vee |\bar{x}|\vee |\bar{y}|\leq R$.
\end{Assu}
\begin{Assu}\label{assumption2}
There exist positive parameter $p\in[2,\infty)$ and positive constant $K_{1}$ such that
\begin{align}\label{a6}
x^{T}f(x,y)+\frac{p-1}{2}|g(x,y)|^{2}\leq K_{1}(1+|x|^{2}+|y|^{2}).
\end{align}
\end{Assu}

\begin{Lem}(\cite{MR2158891})
Assume that Assumption \ref{assumption1} and \ref{assumption2} hold. Then for any given initial data, there is a unique global solution $X(t)$ to \eqref{a1} on $t\in [-\tau, \infty).$
Moreover, the solution has the property that
\begin{align}
\sup_{-\tau\leq t\leq T}\mathbb{E}|X(t)|^{p}<\infty.
\end{align}
\end{Lem}

Next, the concepts of C-stability and B-consistency in \cite{MR3493491} are modified appropriately and the corresponding definitions for SSDEs are given as follows.
\begin{Defi}
A stochastic one-step method $(\Psi, h,\xi)$ for SDDEs \eqref{a1} is said to be stochastic C-stable if for $\eta\in(1,\infty)$ and all random variables $Y, Z\in L^{2}(\Omega,\mathscr{F}_{t},\mathbb{P};\mathbb{R}^{d})$
\begin{align}\label{a2}
\begin{split}
\big{\|}\mathbb{E}&[\Psi(Y,\bar{Y},h)-\Psi(Z,\bar{Z},h)|\mathscr{F}_{t}]\big{\|}_{L^{2}(\Omega;\mathbb{R}^{d})}^{2}\\
&+\eta\big{\|}(id-\mathbb{E}[\cdot|\mathscr{F}_{t}])(\Psi(Y,\bar{Y},h)-\Psi(Z,\bar{Z},h))\big{\|}_{L^{2}(\Omega;\mathbb{R}^{d})}^{2}\\
\leq& (1+C_{stab}h)\big{\|}Y-Z\big{\|}_{L^{2}(\Omega;\mathbb{R}^{d})}^{2}+C_{stab}h\big{\|}\bar{Y}-\bar{Z}\big{\|}_{L^{2}(\Omega;\mathbb{R}^{d})}^{2},
\end{split}
\end{align}
\end{Defi}
where $(id-\mathbb{E}[\cdot|\mathscr{F}_{t}])Y=Y-\mathbb{E}[Y|\mathscr{F}_{t}]$.
\begin{Defi}
A stochastic one-step method $(\Psi, h,\xi)$  for SDDEs \eqref{a1} is said to be stochastic B-consistent of order $\gamma$ if
\begin{align}\label{a3}
\big{\|}\mathbb{E}[X(t+h)-\Psi(X(t),X(t-\tau),h)|\mathscr{F}_{t}]\big{\|}_{L^{2}(\Omega;\mathbb{R}^{d})}\leq C_{cons}h^{\gamma+1},
\end{align}
and
\begin{align}\label{a4}
\big{\|}(id-\mathbb{E}[\cdot|\mathscr{F}_{t}])\big{(}X(t+h)-\Psi(X(t),X(t-\tau),h)\big{)}\big{\|}_{L^{2}(\Omega;\mathbb{R}^{d})}\leq C_{cons}h^{\gamma+\frac{1}{2}}.
\end{align}
\end{Defi}
\section{Convergence Theorem}
The next stability lemma plays an important role in the convergence analysis.
\begin{Lem}
If $(\Psi, h, \xi)$ is stochastically C-stable one-step method with constants $C_{stab}$ and $\eta \in (1,\infty)$,  then for every grid function $Z\in \mathscr{G}^{2}(\mathcal{T}_{h})$,
\begin{align}\label{a7}
\begin{split}
\max_{n\in\{0\cdots N\}}&\|Z(t_{n})-X_{h}(t_{n})\|_{L^{2}(\Omega;\mathbb{R}^{d})}^{2}\\
\leq& e^{2(1+C_{stab}(1+h))T}\bigg{(}\sum_{i=1}^{M}\big{\|}Z(t_{i-M})-\xi(t_{i-M})\big{\|}_{L^{2}(\Omega;\mathbb{R}^{d})}^{2}+
\|Z(t_{0})-X_{h}(t_{0})\|_{L^{2}(\Omega;\mathbb{R}^{d})}^{2}\\
&+\sum_{i=1}^{N}(1+h^{-1})\big{\|}\mathbb{E}\big[Z(t_{i})-\Psi(Z(t_{i-1}),Z(t_{i-M}),h)|\mathscr{F}_{t_{i-1}}]\big{\|}_{L^{2}(\Omega;\mathbb{R}^{d})}^{2}\\
&+C_{\eta}\sum_{i=1}^{N}\big{\|}(id-\mathbb{E}[\cdot|\mathscr{F}_{t_{i-1}}])\big{(}Z(t_{i})-\Psi(Z(t_{i-1}),Z(t_{i-M}),h)\big{)}\big{\|}_{L^{2}(\Omega;\mathbb{R}^{d})}^{2}
 \bigg{)},
\end{split}
\end{align}
where $Z(t_{i-M})$, $\xi(t_{i-M}), i=0, 1, \cdots, M$, are defined by $Z(t_{i}-\tau)$ and $\xi(t_{i}-\tau)$, respectively.
\end{Lem}
\textbf{Proof}. For every $1\leq i\leq N$, let $e_{h}(t_{i}):=Z(t_{i})-X_{h}(t_{i})$,
\begin{align*}
\|e_{h}(t_{i})\|_{L^{2}(\Omega;\mathbb{R}^{d})}^{2}=\big{\|}\mathbb{E}[e_{h}(t_{i})|\mathscr{F}_{t_{i-1}}]\big{\|}_{L^{2}(\Omega;\mathbb{R}^{d})}^{2}
+\big{\|}e_{h}(t_{i})-\mathbb{E}[e_{h}(t_{i})|\mathscr{F}_{t_{i-1}}]\big{\|}_{L^{2}(\Omega;\mathbb{R}^{d})}^{2}
\end{align*}
On account of
\begin{align*}
e_{h}(t_{i})=Z(t_{i})-\Psi(Z(t_{i-1}),Z(t_{i-M}),h)+\Psi(Z(t_{i-1}),Z(t_{i-M}),h)-X_{h}(t_{i}),
\end{align*}
we have
\begin{align*}
\big{\|}\mathbb{E}[e_{h}(t_{i})|\mathscr{F}_{t_{i-1}}]\big{\|}_{L^{2}(\Omega;\mathbb{R}^{d})}\leq &\big{\|}\mathbb{E}[Z(t_{i})-\Psi(Z(t_{i-1}),Z(t_{i-M}),h)|\mathscr{F}_{t_{i-1}}]\big{\|}_{L^{2}(\Omega;\mathbb{R}^{d})}\\
&+\big{\|}\mathbb{E}[\Psi(Z(t_{i-1}),Z(t_{i-M}),h)-X_{h}(t_{i})|\mathscr{F}_{t_{i-1}}]\big{\|}_{L^{2}(\Omega;\mathbb{R}^{d})}.
\end{align*}
By the inequality $(a+b)^{2}=a^2+2ab+b^{2}\leq (1+h^{-1})a^{2}+(1+h)b^{2}$, one may derive that
\begin{align*}
\big{\|}\mathbb{E}[e_{h}(t_{i})|\mathscr{F}_{t_{i-1}}]\big{\|}_{L^{2}(\Omega;\mathbb{R}^{d})}^{2}\leq &(1+h^{-1})\big{\|}\mathbb{E}[Z(t_{i})-\Psi(Z(t_{i-1}),Z(t_{i-M}),h)|\mathscr{F}_{t_{i-1}}]\big{\|}_{L^{2}(\Omega;\mathbb{R}^{d})}^{2}\\
&+(1+h)\big{\|}\mathbb{E}[\Psi(Z(t_{i-1}),Z(t_{i-M}),h)-X_{h}(t_{i})|\mathscr{F}_{t_{i-1}}]\big{\|}_{L^{2}(\Omega;\mathbb{R}^{d})}^{2}.
\end{align*}
Repeating the same process for the item $\big{\|}e_{h}(t_{i})-\mathbb{E}[e_{h}(t_{i})|\mathscr{F}_{t_{i-1}}]\big{\|}_{L^{2}(\Omega;\mathbb{R}^{d})}^{2}$, and replacing $h$ with $\eta-1$, then
\begin{align*}
\big{\|}e_{h}(t_{i})-\mathbb{E}[e_{h}(t_{i})|\mathscr{F}_{t_{i-1}}]\big{\|}_{L^{2}(\Omega;\mathbb{R}^{d})}^{2}\leq &C_{\eta}\big{\|}(id-\mathbb{E}[\cdot|\mathscr{F}_{t_{i-1}}])\big{(}Z(t_{i})-\Psi(Z(t_{i-1}),Z(t_{i-M}),h)\big{)}\big{\|}_{L^{2}(\Omega;\mathbb{R}^{d})}^{2}\\
&+\eta\big{\|}(id-\mathbb{E}[\cdot|\mathscr{F}_{t_{i-1}}])\big{(}\Psi(Z(t_{i-1}),Z(t_{i-M}),h)-X_{h}(t_{i})\big{)}\big{\|}_{L^{2}(\Omega;\mathbb{R}^{d})}^{2},
\end{align*}
where $C_{\eta}=1+(\eta-1)^{-1}$. Consequently, for $1\leq i\leq N$,
\begin{align*}
\|Z(t_{i})&-X_{h}(t_{i})\|_{L^{2}(\Omega;\mathbb{R}^{d})}^{2}\\
\leq& (1+h^{-1})\big{\|}\mathbb{E}[Z(t_{i})-\Psi(Z(t_{i-1}),Z(t_{i-M}),h)|\mathscr{F}_{t_{i-1}}]\big{\|}_{L^{2}(\Omega;\mathbb{R}^{d})}^{2}\\
&+(1+h)\big{\|}\mathbb{E}[\Psi(Z(t_{i-1}),Z(t_{i-M}),h)-X_{h}(t_{i})|\mathscr{F}_{t_{i-1}}]\big{\|}_{L^{2}(\Omega;\mathbb{R}^{d})}^{2}\\
&+C_{\eta}\big{\|}(id-\mathbb{E}[\cdot|\mathscr{F}_{t_{i-1}}])(Z(t_{i})-\Psi(Z(t_{i-1}),Z(t_{i-M}),h))\big{\|}_{L^{2}(\Omega;\mathbb{R}^{d})}^{2}\\
&+\eta\big{\|}(id-\mathbb{E}[\cdot|\mathscr{F}_{t_{i-1}}])(\Psi(Z(t_{i-1}),Z(t_{i-M}),h)-X_{h}(t_{i}))\big{\|}_{L^{2}(\Omega;\mathbb{R}^{d})}^{2}.
\end{align*}
Using the fact that $X_{h}(t_{i})=\Psi(X_{h}(t_{i-1}),X_{h}(t_{i-M}),h)$ and \eqref{a2}, we have
\begin{align*}
\|Z(t_{i})&-X_{h}(t_{i})\|_{L^{2}(\Omega;\mathbb{R}^{d})}^{2}\\
\leq& (1+h^{-1})\big{\|}\mathbb{E}[Z(t_{i})-\Psi(Z(t_{i-1}),Z(t_{i-M}),h)|\mathscr{F}_{t_{i-1}}]\big{\|}_{L^{2}(\Omega;\mathbb{R}^{d})}^{2}\\
&+C_{\eta}\big{\|}(id-\mathbb{E}[\cdot|\mathscr{F}_{t_{i-1}}])(Z(t_{i})-\Psi(Z(t_{i-1}),Z(t_{i-M}),h))\big{\|}_{L^{2}(\Omega;\mathbb{R}^{d})}^{2}\\
&+\big{(}1+(1+C_{stab}(1+h))h\big{)}\big{\|}Z(t_{i-1})-X_{h}(t_{i-1})\big{\|}_{L^{2}(\Omega;\mathbb{R}^{d})}^{2}\\
&+C_{stab}h(1+h)\big{\|}Z(t_{i-M})-X_{h}(t_{i-M})\big{\|}_{L^{2}(\Omega;\mathbb{R}^{d})}^{2},
\end{align*}
where we have used the inequality
\begin{align*}
&h\big{\|}\mathbb{E}[\Psi(Z(t_{i-1}),Z(t_{i-M}),h)-X_{h}(t_{i})|\mathscr{F}_{t_{i-1}}]\big{\|}_{L^{2}(\Omega;\mathbb{R}^{d})}^{2}\\
&\leq h(1+C_{stab}h)\|Z(t_{i-1})-X_{h}(t_{i-1})\|_{L^{2}(\Omega;\mathbb{R}^{d})}^{2}+hC_{stab}h\|Z(t_{i-M})-X_{h}(t_{i-M})\|_{L^{2}(\Omega;\mathbb{R}^{d})}^{2}.
\end{align*}
Choose sufficiently small $h$ such that $C_{stab}h(1+h)<1$. Then,  summing $i$ over $1$ to $n$ yields
\begin{align}\label{a12}
\begin{split}
\|&Z(t_{n})-X_{h}(t_{n})\|_{L^{2}(\Omega;\mathbb{R}^{d})}^{2}-\|Z(t_{0})-X_{h}(t_{0})\|_{L^{2}(\Omega;\mathbb{R}^{d})}^{2}\\
=&\sum_{i=1}^{n}\bigg{(}\|Z(t_{i})-X_{h}(t_{i})\|_{L^{2}(\Omega;\mathbb{R}^{d})}^{2}-\|Z(t_{i-1})-X_{h}(t_{i-1})\|_{L^{2}(\Omega;\mathbb{R}^{d})}^{2}\bigg{)}\\
\leq& \sum_{i=1}^{n}\bigg{(}(1+h^{-1})\big{\|}\mathbb{E}[Z(t_{i})-\Psi(Z(t_{i-1}),Z(t_{i-M}),h)|\mathscr{F}_{t_{i-1}}]\big{\|}_{L^{2}(\Omega;\mathbb{R}^{d})}^{2}\\
&+C_{\eta}\big{\|}(id-\mathbb{E}[\cdot|\mathscr{F}_{t_{i-1}}])(Z(t_{i})-\Psi(Z(t_{i-1}),Z(t_{i-M}),h))\big{\|}_{L^{2}(\Omega;\mathbb{R}^{d})}^{2}\\
&+2(1+C_{stab}(1+h))h\big{\|}Z(t_{i-1})-X_{h}(t_{i-1})\big{\|}_{L^{2}(\Omega;\mathbb{R}^{d})}^{2}\bigg{)}\\
&+\sum_{i=1}^{M}\|Z(t_{i-M})-\xi(t_{i-M})\|_{L^{2}(\Omega;\mathbb{R}^{d})}^{2}.
\end{split}
\end{align}
Finally, the desired assertion follows from \eqref{a12} and the discrete Gronwall inequality. \qed

The following theorem shows that convergence can be derived from stability plus consistency.
\begin{Theo}
If a stochastic one-step method $(\Psi, h, \xi)$ is stochastic C-stability and B-consistent of order $\gamma$, then there exists a constant $C$ such that
\begin{align*}
\max_{n\in \{0,\cdots,N\}}\|X(t_{n})-X_{h}(t_{n})\|_{L^{2}(\Omega;\mathbb{R}^{d})}\leq Ch^{\gamma},
\end{align*}
where X is the exact solution of \eqref{a1} and $X_{h}$ is the grid function corresponding to $(\Psi, h, \xi)$ with time step $h$.
\end{Theo}
\textbf{Proof}. Due to the fact that $X(t_{i-M})=X_{h}(t_{i-M})=\xi(t_{i-M}), i=0, 1, \cdots, M$,  we obtain
\begin{align*}
&\max_{n\in \{0,\cdots,N\}}\|X(t_{n})-X_{h}(t_{n})\|_{L^{2}(\Omega;\mathbb{R}^{d})}^{2}\\
\leq& e^{2(1+C_{stab}(1+h))T}\bigg{(}\sum_{i=1}^{N}(1+h^{-1})\big{\|}\mathbb{E}[Z(t_{i})-
\Psi(Z(t_{i-1}),Z(t_{i-M}),h)|\mathscr{F}_{t_{i-1}}]\big{\|}_{L^{2}(\Omega;\mathbb{R}^{d})}^{2}\\
&+C_{\eta}\sum_{i=1}^{N}\big{\|}(id-\mathbb{E}[\cdot|\mathscr{F}_{t_{i-1}}])(Z(t_{i})-\Psi(Z(t_{i-1}),Z(t_{i-M}),h))\big{\|}_{L^{2}(\Omega;\mathbb{R}^{d})}^{2}\bigg{)}.
\end{align*}
It follows from \eqref{a3} and \eqref{a4} that
\begin{align*}
&\max_{n\in \{0, \cdots, N\}}\|X(t_{n})-X_{h}(t_{n})\|_{L^{2}(\Omega;\mathbb{R}^{d})}^{2}\\
\leq& e^{2(1+C_{stab}(1+h))T}C_{cons}^{2}\sum_{i=1}^{N}\bigg{(}(1+h^{-1})h^{2(\gamma+1)}+C_{\eta}h^{2\gamma+1}\bigg{)}\\
\leq& Ch^{2\gamma}.
\end{align*}
The proof is completed now.\qed

\section{C-stability and B-consistency of the PEM Method}
We follow the notation in \cite{MR3493491}
$$
x^{\circ}:=\min(1,h^{-\alpha}|x|^{-1})x,
$$
and denote
$$
\bar{x}^{\circ}:=\min(1,h^{-\alpha}|\bar{x}|^{-1})\bar{x},
$$
where $x\in \mathbb{R}^{d}$ and step size $h\in (0,1]$.
\begin{Lem}(\cite{MR3493491})\label{lemma4.1}
For every $\alpha\in (0,\infty)$ and $h\in (0,1]$ the mapping $\mathbb{R}^{d}\ni x |\rightarrow x^{\circ}\in \mathbb{R}^{d}$ is globally Lipschitz continous with Lipschitz
constant $1$, i.e.,
$$
|x_{1}^{\circ}-x_{2}^{\circ}|\leq |x_{1}-x_{2}|
$$
for all $x_{1}, x_{2}\in \mathbb{R}^{d}$.
\end{Lem}

\begin{Lem}
If Assumption \ref{assu1} is fulfilled with $L\in(0,\infty)$, $q\in (1,\infty)$ and $\eta\in (\frac{1}{2},\infty)$, then the functions $x^{\circ}$, $\bar{x}^{\circ}$ with parameter $\alpha\in (0,\frac{1}{2(q-1)})$ and $h\in (0,1]$ satisfy
\begin{align*}
\big{|}x_{1}^{\circ}&-x_{2}^{\circ}+h(f(x_{1}^{\circ},\bar{x_{1}}^{\circ})-f(x_{2}^{\circ},\bar{x_{2}}^{\circ}))\big{|}^{2}
+2\eta h\big{|}g(x_{1}^{\circ},\bar{x_{1}}^{\circ})-g(x_{2}^{\circ},\bar{x_{2}}^{\circ})\big{|}^{2}\\
&\leq (1+Ch)|x_{1}-x_{2}|^{2}+Ch|\bar{x_{1}}-\bar{x_{2}}|^{2}
\end{align*}
for all $x_{1}, x_{2}\in \mathbb{R}^{d}$.
\end{Lem}
\textbf{Proof}. By \eqref{a13}, we obtain that
\begin{align*}
\big{|}x_{1}^{\circ}&-x_{2}^{\circ}+h(f(x_{1}^{\circ},\bar{x_{1}}^{\circ})-f(x_{2}^{\circ},\bar{x_{2}}^{\circ}))\big{|}^{2}\\
=&|x_{1}^{\circ}-x_{2}^{\circ}|^{2}+2h\langle x_{1}^{\circ}-x_{2}^{\circ},f(x_{1}^{\circ},\bar{x_{1}}^{\circ})-f(x_{2}^{\circ},\bar{x_{2}}^{\circ})\rangle
+h^{2}|f(x_{1}^{\circ},\bar{x_{1}}^{\circ})-f(x_{2}^{\circ},\bar{x_{2}}^{\circ})|^{2}\\
\leq& (1+2Lh)|x_{1}^{\circ}-x_{2}^{\circ}|^{2}+2Lh|\bar{x_{1}}^{\circ}-\bar{x_{2}}^{\circ}|^{2}-2\eta h|g(x_{1}^{\circ},\bar{x_{1}}^{\circ})
-g(x_{2}^{\circ},\bar{x_{2}}^{\circ})|^{2}\\
&+h^{2}|f(x_{1}^{\circ},\bar{x_{1}}^{\circ})-f(x_{2}^{\circ},\bar{x_{2}}^{\circ})|^{2}.
\end{align*}
Note that
\begin{align*}
&|f(x_{1}^{\circ},\bar{x_{1}}^{\circ})-f(x_{2}^{\circ},\bar{x_{2}}^{\circ})|\\
&\leq L(1+|x_{1}^{\circ}|^{q-1}+|\bar{x_{1}}^{\circ}|^{q-1}+|x_{2}^{\circ}|^{q-1}+|\bar{x_{2}}^{\circ}|^{q-1})(|x_{1}^{\circ}-x_{2}^{\circ}|+|\bar{x_{1}}^{\circ}-\bar{x_{2}}^{\circ}|)\\
&\leq L(1+4h^{-\alpha(q-1)})(|x_{1}-x_{2}|+|\bar{x_{1}}-\bar{x_{2}}|)\\
&\leq L(1+4h^{-1/2})(|x_{1}-x_{2}|+|\bar{x_{1}}-\bar{x_{2}}|),
\end{align*}
where we have used \eqref{a9}, Lemma \ref{lemma4.1}, $|x_{1}^{\circ}|, |x_{2}^{\circ}|, |\bar{x_{1}}^{\circ}|, |\bar{x_{2}}^{\circ}|\leq h^{-\alpha}$ and $\alpha\in (0, \frac{1}{2(q-1)}]$. Consequently,
\begin{align*}
&\big{|}x_{1}^{\circ}-x_{2}^{\circ}+h(f(x_{1}^{\circ},\bar{x_{1}}^{\circ})-f(x_{2}^{\circ},\bar{x_{2}}^{\circ}))\big{|}^{2}
+2\eta h\big{|}g(x_{1}^{\circ},\bar{x_{1}}^{\circ})-g(x_{2}^{\circ},\bar{x_{2}}^{\circ}))\big{|}^{2}\\
&\leq (1+2Lh)|x_{1}^{\circ}-x_{2}^{\circ}|^{2}+2Lh|\bar{x_{1}}^{\circ}-\bar{x_{2}}^{\circ}|^{2}+h^{2} 2L^{2}(1+4h^{-1/2})^{2}(|x_{1}-x_{2}|^{2}+|\bar{x_{1}}-\bar{x_{2}}|^{2})\\
&\leq (1+Ch)|x_{1}-x_{2}|^{2}+Ch|\bar{x_{1}}-\bar{x_{2}}|^{2}.
\end{align*}
The direct application of the above lemma can deduce that the projected Euler method is C-stable.

Next, we show the PEM method is B-consistent of order $1/2$.
\begin{Lem}\label{lemma4.3}
If Assumption \ref{assu1} is fulfilled with $L\in (0,\infty)$, $q\in(1,\infty)$, and $\sup\limits_{t\in[-\tau,T]}\|X(t)\|_{L^{pq}(\Omega,\mathbb{R}^{d})}<\infty$ holds for some positive constant $p\in [2,\infty)$, then
\begin{align*}
\|X(r_{1})-X(r_{2})\|_{L^{p}(\Omega;\mathbb{R}^{d})}\leq C\bigg{(}1+2\sup_{t\in[-\tau,T]}\|X(t)\|_{L^{pq}(\Omega;\mathbb{R}^{d})}^{q}\bigg{)}~|r_{1}-r_{2}|^{1/2},
\end{align*}
for all $r_{1}, r_{2}\in [0,T]$.
\end{Lem}
\textbf{Proof}. The proof can be deduced from Proposition $5.4$ of \cite{MR3493491} easily. In fact, we just need replace $\sup\limits_{t\in[0,T]}\|X(t)\|_{L^{pq}(\Omega;\mathbb{R}^{d})}^{q}$ with $2\sup\limits_{t\in[-\tau,T]}\|X(t)\|_{L^{pq}(\Omega;\mathbb{R}^{d})}^{q}$. So we omit the detail of
proof here. \qed

\begin{Lem}\label{lemma4.4}
If $f$ and $g$ satisfy Assumption \ref{assu1} with $L\in (0,\infty)$ and $q\in (1,\infty)$, the exact solution of \eqref{a1} satisfy $\sup\limits_{t\in[-\tau,T]}\|X(t)\|_{L^{4q-2}(\Omega;\mathbb{R}^{d})}<\infty$, then for $\forall s_{1}\in [r_{1}, r_{2}]$, there exists a constant $C$ such that
\begin{align*}
&\int_{r_{1}}^{r_{2}}\|f(X(s),X(s-\tau))-f(X(s_{1}),X(s_{1}-\tau))\|_{L^{2}(\Omega;\mathbb{R}^{d})}ds\\
&\leq C\bigg{(}1+4\sup_{t\in[-\tau,T]}\|X(t)\|_{L^{4q-2}(\Omega;\mathbb{R}^{d})}^{2q-1}\bigg{)}~|r_{1}-r_{2}|^{3/2},
\end{align*}
for all $r_{1}, r_{2}\in [0,T]$.
\end{Lem}
\textbf{Proof}. By \eqref{a9} and H\"{o}lder's inequality, we get
\begin{align}\label{a10}
\begin{split}
&\|f(X(s),X(s-\tau))-f(X(s_{1}),X(s_{1}-\tau))\|_{L^{2}(\Omega;\mathbb{R}^{d})}\\
&\leq L\big{\|}(1+|X(s)|^{q-1}+|X(s-\tau)|^{q-1}+|X(s_{1})|^{q-1}+|X(s_{1}-\tau)|^{q-1})\big{(}|X(s)-X(s_{1})|\big{)}\big{\|}_{L^{2}(\Omega;\mathbb{R}^{d})}\\
&+ L\big{\|}(1+|X(s)|^{q-1}+|X(s-\tau)|^{q-1}+|X(s_{1})|^{q-1}+|X(s_{1}-\tau)|^{q-1})\big{(}|X(s-\tau)-X(s_{1}-\tau)|\big{)}\big{\|}_{L^{2}(\Omega;\mathbb{R}^{d})}\\
&\leq L\bigg{(}1+4\sup_{t\in[-\tau,T]}\|X(t)\|_{L^{2\rho^{'}(q-1)}(\Omega;\mathbb{R}^{d})}^{q-1}\bigg{)}~\|X(s)-X(s_{1})\|_{L^{2\rho}(\Omega;\mathbb{R}^{d})}\\
&~~~+L\bigg{(}1+4\sup_{t\in[-\tau,T]}\|X(t)\|_{L^{2\rho^{'}(q-1)}(\Omega;\mathbb{R}^{d})}^{q-1}\bigg{)}~\|X(s-\tau)-X(s_{1}-\tau)\|_{L^{2\rho}(\Omega;\mathbb{R}^{d})},
\end{split}
\end{align}
where $\rho=2-\frac{1}{q}$ and $\rho^{'}=\frac{2q-1}{q-1}$.

Without loss of generality, we discuss the last term in \eqref{a10} with three different cases,\\
$Case~ 1$: $s-\tau>0$ and $s_{1}-\tau<0$,
$$
|X(s-\tau)-X(s_{1}-\tau)|\leq |X(s-\tau)-X(0)|+|X(0)-X(s_{1}-\tau)|,
$$
hence,
\begin{align*}
&\|X(s-\tau)-X(s_{1}-\tau)\|_{L^{2\rho}(\Omega;\mathbb{R}^{d})}\\
&\leq \|X(s-\tau)-X(0)\|_{L^{2\rho}(\Omega;\mathbb{R}^{d})}+\|\xi(0)-\xi(s_{1}-\tau)\|_{L^{2\rho}(\Omega;\mathbb{R}^{d})}\\
&\leq C\bigg{(}1+2\sup_{t\in[-\tau,T]}\|X(t)\|_{L^{2\rho q}(\Omega;\mathbb{R}^{d})}^{q}\bigg{)}~|s-\tau|^{1/2}+K_{1}|s_{1}-\tau|^{\beta}\\
&\leq C\bigg{(}1+2\sup_{t\in[-\tau,T]}\|X(t)\|_{L^{4q-2}(\Omega;\mathbb{R}^{d})}^{q}\bigg{)}~|s-s_{1}|^{1/2}+K_{1}|s_{1}-s|^{\beta}\\
&\leq C\bigg{(}1+2\sup_{t\in[-\tau,T]}\|X(t)\|_{L^{4q-2}(\Omega;\mathbb{R}^{d})}^{q}\bigg{)}~|r_{1}-r_{2}|^{1/2}+K_{1}|r_{1}-r_{2}|^{\beta},
\end{align*}
where Assumption \ref{assu2} and Lemma \ref{lemma4.3} are used.\\
$Case~ 2$: $s-\tau>0$ and $s_{1}-\tau>0$, it follows from Lemma \ref{lemma4.3} that
\begin{align*}
&\|X(s-\tau)-X(s_{1}-\tau)\|_{L^{2\rho}(\Omega;\mathbb{R}^{d})}\\
&\leq C\bigg{(}1+2\sup_{t\in[-\tau,T]}\|X(t)\|_{L^{2\rho q}(\Omega;\mathbb{R}^{d})}^{q}\bigg{)}~|(s-\tau)-(s_{1}-\tau)|^{1/2}\\
&\leq C\bigg{(}1+2\sup_{t\in[-\tau,T]}\|X(t)\|_{L^{4q-2}(\Omega;\mathbb{R}^{d})}^{q}\bigg{)}~|s-s_{1}|^{1/2}\\
&\leq C\bigg{(}1+2\sup_{t\in[-\tau,T]}\|X(t)\|_{L^{4q-2}(\Omega;\mathbb{R}^{d})}^{q}\bigg{)}~|r_{2}-r_{1}|^{1/2}.
\end{align*}
$Case~ 3$: $s-\tau<0$ and $s_{1}-\tau<0$,  Assumption \ref{assu2} implies that
\begin{align*}
&\|X(s-\tau)-X(s_{1}-\tau)\|_{L^{2\rho}(\Omega;\mathbb{R}^{d})}\\
&\leq \|\xi(s-\tau)-\xi(s_{1}-\tau)\|_{L^{2\rho}(\Omega;\mathbb{R}^{d})}\\
&\leq K_{1}|(s-\tau)-(s_{1}-\tau)|^{\beta}\leq K_{1}|r_{2}-r_{1}|^{\beta}.
\end{align*}
Together \eqref{a10} with three cases above, we have
\begin{align*}
\|f(X(s),X(s-\tau))-f(X(s_{1}),X(s_{1}-\tau))\|_{L^{2}(\Omega;\mathbb{R}^{d})}\leq C\bigg{(}1+4\sup_{t\in[-\tau,T]}\|X(t)\|_{L^{4q-2}(\Omega;\mathbb{R}^{d})}^{2q-1}\bigg{)}~|r_{1}-r_{2}|^{1/2}.
\end{align*}
So the desired assertion follows. \qed

\begin{Lem}\label{lemma4.5}
If the coefficients $f$ and $g$ satisfy Assumption \ref{assu1} with $L\in (0,\infty)$ and $q\in (1,\infty)$, the exact solution of \eqref{a1} satisfy $\sup\limits_{t\in[-\tau,T]}\|X(t)\|_{L^{4q-2}(\Omega;\mathbb{R}^{d})}<\infty$, then there exists a constant $C$ such that
\begin{align}\label{a11}
\begin{split}
&\big{\|}\int_{r_{1}}^{r_{2}}g(X(s),X(s-\tau))-g(X(r_{1}),X(r_{1}-\tau))dW(s)\big{\|}_{L^{2}(\Omega;\mathbb{R}^{d})}\\
&\leq C\bigg{(}1+4\sup_{t\in[-\tau,T]}\|X(t)\|_{L^{4q-2}(\Omega;\mathbb{R}^{d})}^{2q-1}\bigg{)}~|r_{1}-r_{2}|,
\end{split}
\end{align}
for all $r_{1}, r_{2}\in [0,T]$.
\end{Lem}
\textbf{Proof}. It\^{o} isometry formula yields
\begin{align*}
&\big{\|}\int_{r_{1}}^{r_{2}}g(X(s),X(s-\tau))-g(X(r_{1}),X(r_{1}-\tau))dW(s)\big{\|}_{L^{2}(\Omega;\mathbb{R}^{d})}\\
&\leq \bigg{(}\int_{r_{1}}^{r_{2}}\|g(X(s),X(s-\tau))-g(X(r_{1}),X(r_{1}-\tau))\|_{L^{2}(\Omega;\mathbb{R}^{d})}^{2}ds\bigg{)}^{1/2}.
\end{align*}
Repeating the proof in Lemma \ref{lemma4.4}, and we get
\begin{align*}
\begin{split}
&\|g(X(s),X(s-\tau))-g(X(r_{1}),X(r_{1}-\tau))\|_{L^{2}(\Omega;\mathbb{R}^{d})}\\
&\leq L\big{\|}(1+|X(s)|^{q-1}+|X(s-\tau)|^{q-1}+|X(r_{1})|^{q-1}+|X(r_{1}-\tau)|^{q-1})\big{(}|X(s)-X(r_{1})|\big{)}\big{\|}_{L^{2}(\Omega;\mathbb{R}^{d})}\\
&~~~+ L\big{\|}(1+|X(s)|^{q-1}+|X(s-\tau)|^{q-1}+|X(r_{1})|^{q-1}+|X(r_{1}-\tau)|^{q-1})\big{(}|X(s-\tau)-X(r_{1}-\tau)|\big{)}\big{\|}_{L^{2}(\Omega;\mathbb{R}^{d})}\\
&\leq L\bigg{(}1+4\sup_{t\in[-\tau,T]}\|X(t)\|_{L^{2\rho'(q-1)}(\Omega;\mathbb{R}^{d})}^{q-1}\bigg{)}~\|X(s)-X(r_{1})\|_{L^{2\rho}(\Omega;\mathbb{R}^{d})}\\
&~~~+L\bigg{(}1+4\sup_{t\in[-\tau,T]}\|X(t)\|_{L^{2\rho'(q-1)}(\Omega;\mathbb{R}^{d})}^{q-1}\bigg{)}~\|X(s-\tau)-X(r_{1}-\tau)\|_{L^{2\rho}(\Omega;\mathbb{R}^{d})}\\
&\leq C\bigg{(}1+4\sup_{t\in[-\tau,T]}\|X(t)\|_{L^{4q-2}(\Omega;\mathbb{R}^{d})}^{2q-1}\bigg{)}~|r_{1}-r_{2}|^{1/2}.
\end{split}
\end{align*}
The required inequality \eqref{a11} follows. \qed

Before giving the consistency results, we first present the following key lemma similar to Lemma $6.5$ of \cite{MR3493491}.
\begin{Lem}\label{lemma 4.6}
Denote $L\in (0,\infty)$ and $\kappa\in [1,\infty)$. If for $p\in (2,\infty)$, $Y\in L^{p\kappa}(\Omega;\mathbb{R}^{d})$, and the measurable mapping $\varphi: \mathbb{R}^{d}\rightarrow \mathbb{R}^{d}$ has  the following properties
$$
|\varphi(x,\bar{x})|\leq L(1+|x|^{\kappa}+|\bar{x}|^{\kappa}),
$$
then there exists a constant $C$ which depends on $p, L$, but not on $h$ such that for all $h\in (0,1]$
$$
\|\varphi(Y,\bar{Y})-\varphi(Y^{\circ},\bar{Y}^{\circ})\|_{L^{2}(\Omega;\mathbb{R}^{d})}\leq C\big{(}1+\|Y\|_{L^{p\kappa}(\Omega;\mathbb{R}^{d})}^{\kappa}+\|\bar{Y}\|_{L^{p\kappa}(\Omega;\mathbb{R}^{d})}^{\kappa}\big{)}^{p/2}h^{\frac{1}{2}\alpha(p-2)\kappa},
$$
where $\alpha\in (0,\infty)$, $\bar{Y}^{\circ}=\min(1,h^{-\alpha}|\bar{Y}|^{-1})\bar{Y}$ and $Y^{\circ}=\min(1,h^{-\alpha}|Y|^{-1})Y$.
\end{Lem}
\textbf{Proof}. Denote the following measurable sets by
\begin{align*}
A_{h}:=\{\omega\in\Omega:|Y(\omega)|\leq h^{-\alpha}\}\in \mathscr{F},\\
\bar{A}_{h}:=\{\omega\in\Omega:|\bar{Y}(\omega)|\leq h^{-\alpha}\}\in \mathscr{F}.
\end{align*}
Let $B_{h}:=A_{h}\cap\bar{A}_{h}$, and $B_{h}^{c}:=\Omega\setminus B_{h}$, then one can see that
$$
\|\varphi(Y,\bar{Y})-\varphi(Y^{\circ},\bar{Y}^{\circ})\|_{L^{2}(\Omega;\mathbb{R}^{d})}^{2}=\int_{\Omega}|\varphi(Y,\bar{Y})-\varphi(Y^{\circ},\bar{Y}^{\circ})|^2
\mathds{1}_{B_{h}^{c}}(\omega)d\mathbb{P}(\omega).
$$
By the Young inequality $ab\leq \frac{h^{v}}{\rho}a^{\rho}+\frac{1}{\rho'}h^{-v\frac{\rho'}{\rho}}b^{\rho'}$ with $v, \rho=\frac{p}{2}, \rho'=\frac{p}{p-2}\in (0,\infty)$, we have
\begin{align*}
&\int_{\Omega}|\varphi(Y,\bar{Y})-\varphi(Y^{\circ},\bar{Y}^{\circ})|^2\mathds{1}_{B_{h}^{c}}(\omega)d\mathbb{P}(\omega).\\
&\leq \frac{2h^{v}}{p}\|\varphi(Y,\bar{Y})-\varphi(Y^{\circ},\bar{Y}^{\circ})\|_{L^{p}(\Omega;\mathbb{R}^{d})}^{p}+\bigg{(}1-\frac{2}{p}\bigg{)}h^{-\frac{2v}{p-2}}
\mathbb{P}(B_{h}^{c}),
\end{align*}
Furthermore,
\begin{align*}
&\|\varphi(Y,\bar{Y})-\varphi(Y^{\circ},\bar{Y}^{\circ})\|_{L^{p}(\Omega;\mathbb{R}^{d})}\\
&\leq \|\varphi(Y,\bar{Y})\|_{L^{p}(\Omega;\mathbb{R}^{d})}+\|\varphi(Y^{\circ},\bar{Y}^{\circ})\|_{L^{p}(\Omega;\mathbb{R}^{d})}\\
&\leq 2L(1+\|Y\|_{L^{p\kappa}(\Omega;\mathbb{R}^{d})}^{\kappa}+\|\bar{Y}\|_{L^{p\kappa}(\Omega;\mathbb{R}^{d})}^{\kappa}).
\end{align*}
Besides,
\begin{align*}
\mathbb{P}(B_{h}^{c})&=\mathbb{E}[\mathds{1}_{B_{h}^{c}}]\\
&\leq h^{\alpha p\kappa}\mathbb{E}[\mathds{1}_{B_{h}^{c}}|Y|^{p\kappa}]+h^{\alpha p\kappa}\mathbb{E}[\mathds{1}_{B_{h}^{c}}|\bar{Y}|^{p\kappa}]\\
&\leq h^{\alpha p\kappa}(\|Y\|_{L^{p\kappa}(\Omega;\mathbb{R}^{d})}^{p\kappa}+\|\bar{Y}\|_{L^{p\kappa}(\Omega;\mathbb{R}^{d})}^{p\kappa}).
\end{align*}
Let $\alpha p\kappa-\frac{2v}{p-2}=v$, i.e., $v=\alpha(p-2)\kappa$, we obtain
\begin{align*}
\|\varphi(Y,\bar{Y})-\varphi(Y^{\circ},\bar{Y}^{\circ})\|_{L^{2}(\Omega;\mathbb{R}^{d})}^{2}\leq& \frac{2}{p}(2L)^{p}h^{\alpha(p-2)\kappa}(1+\|Y\|_{L^{p\kappa}(\Omega;\mathbb{R}^{d})}^{\kappa}+\|\bar{Y}\|_{L^{p\kappa}(\Omega;\mathbb{R}^{d})}^{\kappa})^{p}\\
&+ \bigg{(}1-\frac{2}{p}\bigg{)}h^{\alpha(p-2)\kappa}(\|Y\|_{L^{p\kappa}(\Omega;\mathbb{R}^{d})}^{p\kappa}+\|\bar{Y}\|_{L^{p\kappa}(\Omega;\mathbb{R}^{d})}^{p\kappa}),
\end{align*}
which completes the proof.\qed

We conclude this section with a theorem of B-consistency of the projected Euler method.
\begin{Theo}
If the coefficients $f$ and $g$ satisfy Assumption \ref{assu1} with $L\in (0,\infty)$ and $q\in (1,\infty)$, the exact solution of \eqref{a1} satisfy $\sup\limits_{t\in[-\tau,T]}\|X(t)\|_{L^{6q-4}(\Omega;\mathbb{R}^{d})}<\infty$, then the order of B-consistent of the projected Euler method $(\Psi^{PEM}, h, \xi)$ with
$\alpha=\frac{1}{2(q-1)}$  is $\frac{1}{2}$.
\end{Theo}
\textbf{Proof}. By \eqref{a1}, \eqref{a21}, we have
\begin{align*}
X(t+h)-\Psi^{PEM}(X(t),X(t-h),h)=&\int_{t}^{t+h}f(X(s),X(s-\tau))-f(X(t),X(t-\tau))ds\\
&+X(t)+h f(X(t),X(t-\tau))-X^{\circ}(t)-h f(X^{\circ}(t),X^{\circ}(t-\tau))\\
&+\int_{t}^{t+h} g(X(s),X(s-\tau))-g(X(t),X(t-\tau))dW(s)\\
&+\big{(}g(X(t),X(t-\tau))-g(X^{\circ}(t),X^{\circ}(t-\tau))\big{)}\Delta_{h}W(t),
\end{align*}
where $X^{\circ}(t)=\min\big{(}1,h^{-\alpha}|X(t)|^{-1}\big{)}X(t)$, $X^{\circ}(t-\tau)=\min\big{(}1,h^{-\alpha}|X(t-\tau)|^{-1}\big{)}X(t-\tau)$. Moreover,
\begin{align}\label{a14}
\begin{split}
\|\mathbb{E}&[X(t+h)-\Psi^{PEM}(X(t), X(t-h), h)|\mathscr{F}_{t}]\|_{L^{2}(\Omega;\mathbb{R}^{d})}\\
\leq& \int_{t}^{t+h} \big{\|}\mathbb{E}[f(X(s),X(s-\tau))-f(X(t),X(t-\tau))|\mathscr{F}_{t}]\big{\|}_{L^{2}(\Omega;\mathbb{R}^{d})}ds\\
&+\|X(t)-X^{\circ}(t)\|_{L^{2}(\Omega;\mathbb{R}^{d})}+h\|f(X(t),X(t-\tau))-f(X^{\circ}(t),X^{\circ}(t-\tau))\|_{L^{2}(\Omega;\mathbb{R}^{d})}
\end{split}
\end{align}
From the inequality $\|\mathbb{E}[Y|\mathscr{F}_{t}]\|_{L^{2}(\Omega;\mathbb{R}^{d})}\leq \|Y\|_{L^{2}(\Omega;\mathbb{R}^{d})}$ and Lemma \ref{lemma4.4}, we have
$$
\int_{t}^{t+h} \big{\|}\mathbb{E}[f(X(s),X(s-\tau))-f(X(t),X(t-\tau))|\mathscr{F}_{t}]\big{\|}_{L^{2}(\Omega;\mathbb{R}^{d})}ds\leq C_{cons}h^{\frac{3}{2}}.
$$
Applying Lemma \ref{lemma 4.6} to the term $\|X(t)-X^{\circ}(t)\|_{L^{2}(\Omega;\mathbb{R}^{d})}$ with $\kappa=1$ and $p=6q-4$ yields
$$
\|X(t)-X^{\circ}(t)\|_{L^{2}(\Omega;\mathbb{R}^{d})}\leq (1+\|X(t)\|_{L^{6q-4}(\Omega;\mathbb{R}^{d})}+\|X(t-\tau)\|_{L^{6q-4}(\Omega;\mathbb{R}^{d})})^{3q-2}h^{\frac{3}{2}}.
$$
An application of Lemma \ref{lemma 4.6} to $\|f(X(t),X(t-\tau))-f(X^{\circ}(t),X^{\circ}(t-\tau))\|_{L^{2}(\Omega;\mathbb{R}^{d})}$ with $\kappa=q$ and $p=4-\frac{2}{q}$ yields
$$
\|f(X(t),X(t-\tau))-f(X^{\circ}(t),X^{\circ}(t-\tau))\|_{L^{2}(\Omega;\mathbb{R}^{d})}\leq
(1+\|X(t)\|_{L^{4q-2}(\Omega;\mathbb{R}^{d})}^{q}+\|X(t-\tau)\|_{L^{4q-2}(\Omega;\mathbb{R}^{d})}^{q})^{2-\frac{1}{q}}h^{\frac{1}{2}}.
$$
Altogether, and combined with \eqref{a14}, we can find \eqref{a3} is satisfied with $\gamma=\frac{1}{2}$.

Now, we consider another estimation
\begin{align*}
\|(id&-\mathbb{E}[\cdot|\mathscr{F}_{t}])(X(t+h)-\Psi^{PEM}(X(t),X(t-h),h))\|_{L^{2}(\Omega;\mathbb{R}^{d})}\\
\leq& \int_{t}^{t+h}\|(id-\mathbb{E}[\cdot|\mathscr{F}_{t}])(f(X(s),X(s-\tau))-f(X(t),X(t-\tau)))\|_{L^{2}(\Omega;\mathbb{R}^{d})}ds\\
&+\big{\|}\int_{t}^{t+h} g(X(s),X(s-\tau))-g(X(t),X(t-\tau))dW(s)\big{\|}_{L^{2}(\Omega;\mathbb{R}^{d})}\\
&+\big{\|}\big{(}g(X(t),X(t-\tau))-g(X^{\circ}(t),X^{\circ}(t-\tau))\big{)}\Delta_{h}W(t)\big{\|}_{L^{2}(\Omega;\mathbb{R}^{d})},
\end{align*}
Using the inequality $\|(id-\mathbb{E}[\cdot|\mathscr{F}_{t}])Y\|_{L^{2}(\Omega;\mathbb{R}^{d})}\leq \|Y\|_{L^{2}(\Omega;\mathbb{R}^{d})}$, and Lemma \ref{lemma4.4}, it is easy to see that
$$
\int_{t}^{t+h}\|(id-\mathbb{E}[\cdot|\mathscr{F}_{t}])(f(X(s),X(s-\tau))-f(X(t),X(t-\tau)))\|_{L^{2}(\Omega;\mathbb{R}^{d})}ds\leq C_{cons}h^{\frac{3}{2}}.
$$
It follows from Lemma \ref{lemma4.5} that
$$
\big{\|}\int_{t}^{t+h} g(X(s),X(s-\tau))-g(X(t),X(t-\tau))~dW(s)\big{\|}_{L^{2}(\Omega;\mathbb{R}^{d})}\leq C_{cons}h.
$$
In addition,
\begin{align*}
&\big{\|}\big{(}g(X(t),X(t-\tau))-g(X^{\circ}(t),X^{\circ}(t-\tau))\big{)}\Delta_{h}W(t)\big{\|}_{L^{2}(\Omega;\mathbb{R}^{d})}^{2}\\
&=h \big{\|}g(X(t),X(t-\tau))-g(X^{\circ}(t),X^{\circ}(t-\tau))\big{\|}_{L^{2}(\Omega;\mathbb{R}^{d})}^{2}.
\end{align*}
And once again, we apply Lemma \ref{lemma 4.6} to $\big{\|}g(X(t),X(t-\tau))-g(X^{\circ}(t),X^{\circ}(t-\tau))\big{\|}_{L^{2}(\Omega;\mathbb{R}^{d})}$ with $\kappa=q$ and $p=4-\frac{2}{q}$, then
$$
\|g(X(t),X(t-\tau))-g(X^{\circ}(t),X^{\circ}(t-\tau))\|_{L^{2}(\Omega;\mathbb{R}^{d})}\leq
\bigg{(}1+\|X(t)\|_{L^{4q-2}(\Omega;\mathbb{R}^{d})}^{q}+\|X(t-\tau)\|_{L^{4q-2}(\Omega;\mathbb{R}^{d})}^{q}\bigg{)}^{2-\frac{1}{q}}h^{\frac{1}{2}}.
$$
In summary, \eqref{a4} holds with $\gamma=\frac{1}{2}$. \qed

\section{Numerical Experiments}
\begin{example}\label{example1}
We consider the following example \cite{MR3758637}
\begin{equation*}
dy(t)=[-2y(t)+y(t-1)-y^{5}(t)]dt+y^{2}(t)dW(t),
\end{equation*}
for $t\geq 0$ with initial data $y(t)=cos(t)$.
\end{example}

\begin{example}\label{example2}
Next, we consider more general SDDE
\begin{equation*}
dy(t)=[-2y(t)+y(t-1)-y^{5}(t)-y^{5}(t-1)]dt+[y^{2}(t)+y^{2}(t-1)]dW(t),
\end{equation*}
for $t\geq 0$ with initial data $y(t)=cos(t)$.
\end{example}

\begin{figure}[!ht]
\centering
\includegraphics[height=8cm,width=12cm]{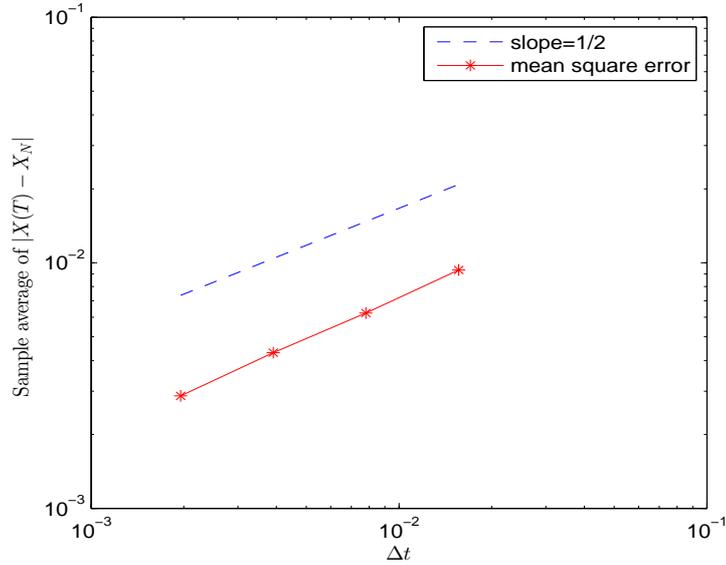}
\caption{Strong convergence of the projected Euler method for Example \ref{example1}.}
\label{1}
\end{figure}

\begin{figure}[!ht]
\centering
\includegraphics[height=8cm,width=12cm]{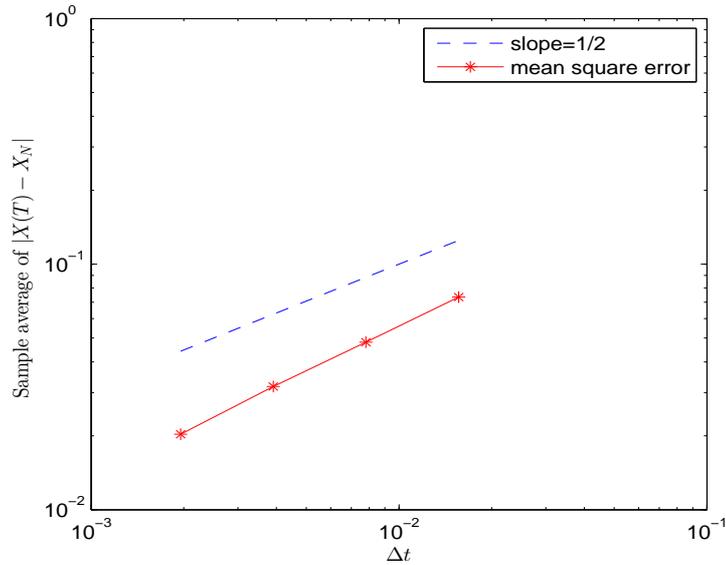}
\caption{Strong convergence of the projected Euler method for Example \ref{example2}.}
\label{2}
\end{figure}

One can see that Example \ref{example1}, \ref{example2} satisfy the condition \eqref{a13}, Assumption \ref{assu1} with $q=5$ and Assumption \ref{assumption2} with $p=30$. Hence, we can choose the projected parameter $\alpha=\frac{1}{2(5-1)}=\frac{1}{8}$. We use discretized Brownian paths over $[0,2]$ with $\Delta t=2^{-13}$. For the sake of simplicity, we regard the projected method with $h=\Delta t$ as good approximation of the exact solution. And compare it with corresponding numerical solution using $h=128\Delta t$, $h=64\Delta t$, $h=32\Delta t$, and $h=16\Delta t$ over $M=1000$ sample paths. We measure the means of absolute errors
at the endpoint $t=T=2$, and denote
\begin{align*}
e_{\Delta t}^{strong}:=\frac{1}{M}\sum_{i=1}^{M}|X_{N}^{(i)}-X(t_{N})^{(i)}|,~~where~~T-\Delta<t_{N}=N\Delta \leq T
\end{align*}
by the endpoint error in the strong sense of the projected Euler method. In Figs \ref{1}, \ref{2}, we plot means of absolute errors $e_{\Delta t}^{strong}$ against $\Delta t$ on log-log scale.
For reference, a dashed blue line is added. We can observe that the convergence rate of the projected Euler method is
$1/2$, which is in accordance with our theoretical results.
\section{Conclusion}
\bibliographystyle{elsarticle-num}
\bibliography{reference}
\end{document}